\documentclass[10pt,twosided]{article}
\usepackage[tbtags]{amsmath}
\usepackage{epsfig,amstext,amssymb,amsthm,latexsym}
\allowdisplaybreaks[4]
\usepackage{color}

\def\dE#1{\textcolor{c30}{#1}}
\def\dE#1{#1}
\def\bE#1{\textcolor{c30}{#1}}
\def\bE#1{#1}

\def\eH#1{\textcolor{c20}{#1}}
\def\eH#1{#1}

\usepackage{amssymb,color}

\definecolor{c20}{rgb}{0.,0.7,0.}
\definecolor{c30}{rgb}{0.,0.,1.}
\definecolor{c40}{rgb}{1,0.1,0.7}
\definecolor{c50}{rgb}{1,0,0}
\definecolor{c60}{rgb}{1,0.9,0.1}

\def\cW#1{{\textcolor{c40}{#1}}}
\def\cW#1{#1}
\def\cZ#1{{\textcolor{c40}{#1}}}
\def\cZ#1{#1}
\def\CC#1{{\textcolor{c40}{#1}}}
\def\CC#1{#1}
\def\CCC#1{{\textcolor{c40}{#1}}}
\def\CCC#1{#1}
\def\ccc#1{{\textcolor{c40}{#1}}}
\def\ccc#1{#1}
\def\dw#1{{\textcolor{c40}{#1}}}
\def\dw#1{#1}
\def\ddw#1{{\textcolor{c40}{#1}}}
\def\ddw#1{#1}

\def\cH#1{\textcolor{c20}{#1}}
\def\cH#1{#1}

\usepackage{amssymb,color}

\newcommand{\ve}{\varepsilon}

\newcommand{\ABs}[1]{ \biggl \lvert #1 \biggr \rvert}

\newcommand{\EE}[1]{\mathbb{E}\left(#1\right)}

\newcommand{\pk}[1]{\mathbb{P} \left( #1 \right) }

\newcommand{\EXP}[1]{\exp \left( #1 \right) }

\newcommand{\R}{\!I\!\!R}

\newcommand{\inr}{\in \R}

\newcommand{\limit}[1]{\lim_{#1 \to   \infty}}

\newcommand{\BQN}{\begin{eqnarray}}
\newcommand{\EQN}{\end{eqnarray}}
\newcommand{\BQNY}{\begin{eqnarray*}}
\newcommand{\EQNY}{\end{eqnarray*}}

\newcommand{\BS}{\begin{sat}}
\newcommand{\ES}{\end{sat}}
\newcommand{\BT}{\begin{theo}}
\newcommand{\ET}{\end{theo}}
\newcommand{\BK}{\begin{korr}}
\newcommand{\EK}{\end{korr}}

\newcommand{\BD}{\begin{de}}
\newcommand{\ED}{\end{de}}

\newcommand{\BRM}{\begin{remarks}}
\newcommand{\ERM}{\end{remarks}}

\newcommand{\BEL}{\begin{lem}}
\newcommand{\EEL}{\end{lem}}

\numberwithin{equation}{section}
\newtheorem{theo}{Theorem}[section]
\newtheorem{sat}[theo]{Proposition}
\newtheorem{de}[theo]{Definition}
\newtheorem{lem}{Lemma}[section]

\newtheorem{korr}[theo]{Corollary}

\newtheorem{remarks}[theo]{Remarks}

\newcommand{\prooftheo}[1]{ \textsc{Proof of Theorem} \ref{#1} }

\newcommand{\prooflem}[1]{\textsc{Proof of Lemma} \ref{#1}}

\newcommand{\QED}{\hfill $\Box$}

\def\d{d}

\def\Cov{\operatorname*{Cov}}
\def\E{\operatorname*{\mathbb{E}}}
\def\I{\operatorname*{\mathbb{I}}}

\newcommand{\netheo}[1]{{Theorem \ref{#1}}}

\newcommand{\nwc}{\newcommand}
\nwc{\COM}[1]{}

\def\IF{\infty}

\def\OH{\overline{H}}

\begin{document}

\centerline{\bf \Large Limit Laws for Maxima of Contracted Stationary Gaussian Sequences}
  \vskip .1 cm

\centerline{Enkelejd Hashorva\footnote{Department of Actuarial Science, Faculty of Business and Economics, University of Lausanne, UNIL-Dorigny 1015 Lausanne, Switzerland, enkelejd.hashorva@unil.ch}  and Zhichao Weng\footnote{Department of Actuarial Science, Faculty of Business and Economics, University of Lausanne, UNIL-Dorigny 1015 Lausanne, Switzerland} }

  \vskip .1 cm
\centerline{Department of Actuarial Science, University of Lausanne}

  \vskip .2 cm

{\bf Abstract}: The principal results of this contribution are the weak and strong limits of maxima of contracted stationary Gaussian random sequences. Due to the random contraction we introduce a modified Berman condition which is sufficient for
the weak convergence of the maxima of the scaled sample. Under a stronger assumption the weak convergence is strengthened to almost convergence.

  \vskip .1 cm


\section{Introduction and Main Result}
If $X,X_n,n\ge 1$ are independent $N(0,1)$ random variables, then it is well-known (see e.g., Berman (1992), Piterbarg (1996) or Falk et al.\ (2010)) that the distribution of sample maxima $M_n=\max_{1 \le i \le n} X_i$
converges (after normalisation) to the Gumbel distribution $\Lambda(x)=\exp(-\exp(-x)),x\inr$, i.e.,
\BQN\label{Mn}
\limit{n}\sup_{x\inr} \ABs{\pk{M_n \le a_n x+ b_n} - \Lambda(x)}=0,
\EQN
where
\BQNY
a_n=(2 \ln   n)^{-\frac{1}{2}}\, \mbox{   and   }\, b_n=(2 \ln   n)^{\frac{1}{2}}-\frac{1}{2}(2 \ln   n)^{-\frac{1}{2}}( \ln   \ln   n + \ln   4 \pi).
\EQNY
Due to some underlying random scaling phenomena, often in applications $Y_i= S_i X_i, i\le n$ are available and not the original observations $X_i,i\le n$, where $S_i$ is some  random factor.
Consider in the following $S, S_n,n\ge 1$ independent non-negative random variables with common distribution function $F$ being independent of $X,X_n,n\ge 1$. We are interested in this paper \dE{in} contraction-type random scaling, i.e., $F$ has a finite upper endpoint, which for simplicity is assumed to be equal to 1.\\
If $S$ is regularly varying at 1 with index $\gamma\ge 0$, i.e.,
 \BQN\label{regvar}
 \limit{u} \frac{\pk{S> 1- t/u}}{\pk{S> 1- 1/u}}=t^\gamma, \quad  \forall t>0,
 \EQN
 then in view of Theorem 3.1 in Hashorva et al.\ (2010) (see also Theorem 4.1 in Hashorva (2013)), the limit relation
 \eqref{Mn} still holds  for  $M_n^*= \max_{1 \le i \le n } S_i X_i$  with constants
 \BQN \label{anbbn}
 b_n= G^{-1}(1- 1/n), \quad a_n=1/b_n \sim (2 \ln n)^{-1/2},
 \EQN
 where $G^{-1}$ is the inverse of the distribution function $G$ of $SX$ and $\sim$ means asymptotical equivalence when $n\to \IF$.
 Our first motivating result states that for any $S$ not equal to 0 the approximation \eqref{Mn} holds.

\BT \label{th0} If $S X$ has distribution function $G$ with generalised inverse $G^{-1}$, then \eqref{Mn} holds for
 $M_n^*$ with constants $a_n,b_n$ as in \eqref{anbbn}.
\ET
The seminal result of Berman \CC{(1964)} shows that if $X_n,n\ge 1$ is a stationary Gaussian sequence with $\rho(n)=\EE{X_1X_n}$, \bE{and $X_1$ is a $N(0,1)$ random variable},
 then the sample maxima $M_n$ still satisfies \eqref{Mn}, provided that the Berman condition
\BQN \label{BMC}
\limit{n} \rho(n) \ln n=0
\EQN
is satisfied. In the sequel we refer to $X_n,n\ge 1$ as a standard stationary Gaussian sequence (ssGs).\\
The main result of this contribution stated below shows that \netheo{th0} can be \bE{stated} for any ssGs, provided that
the Berman condition is accordingly modified, and further some additional restrictions on the random scaling sequence are imposed via the following constrain: \\
{\bf Assumption A}. {\it Let $S$ be a non-negative random variable with distribution function $F$ \CC{which has} upper endpoint 1. \CC{For} any $u \in (\nu, 1)$ with $\nu \in (0,1)$  \BQN\label{eq:bound S}
\cZ{ \pk{S_{\tau}> u}\geq \pk{S> u} \ge  \pk{S_{\gamma}> u}}
\EQN
holds with  $S_{\gamma},S_{\tau}$ two non-negative random variables which have a regularly varying survival function at 1 with non-negative index ${\gamma}$ and ${\tau}$, respectively.
}

We state now our main result:

\BT  \label{th2.1}  If $S$ is such that Assumption A is \bE{satisfied}, then \netheo{th0} holds for any ssGs
$X_n,n\ge 1$ such that for some  $\Delta>2(\gamma-\tau)$
\begin{eqnarray}\label{eq2.1}
\lim_{n \to \infty}\rho(n)( \ln   n)^{1+\Delta}=0.
\end{eqnarray}
\ET

\ddw{This paper is organized as follows}: we continue below with a new Section discussing our main findings and then presenting an extension which strengthens the distributional convergence of maxima
\ccc{$(M_n^* - b_n)/a_n$}
to almost sure convergence. Proofs and auxiliary results are displayed in Section 3.

\section{Discussion and Extensions}
In the light of extreme value theory (see e.g., Resnick (1987), Embrechts  et al.\  (1997), Falk et al.\ (2010)) \CC{the result \eqref{Mn}}  
means that the distribution function $\Phi$ is in the Gumbel max-domain of attraction (MDA). A general univariate distribution function $\CC{H}$  with upper endpoint $\IF$ is in the Gumbel MDA (abbreviated  $H\in GMDA(w))$ if (set $\overline{H}=1-H$)
\BQN\label{G}
\frac{\overline{H}(u+  x/ w(u))}{\overline{H}(u)}\sim \exp(-x), \quad \forall x\ge 0,
\EQN
with $w(\cdot)$ some positive scaling function. Again we write  $\sim$ to mean asymptotic  equivalence of two functions when the argument (typically $u$) approaches infinity. For the standard Gaussian distribution function $\Phi$ on $\R$ we have $\Phi \in GMDA(w)$ where $w(x)=x$.
Consequently, \netheo{th0} means that $SX$ has distribution function $G\in GMDA(w)$ with scaling function $w(x)=x$ \dE{whenever the random variable $S\ge 0$ is  bounded and independent of $X$ which has distribution function $\Phi$.}

Regarding Assumption A we mention that it is satisfied by a large class of random contraction $S$, for instance if
 $S$ is a Beta random variable, or $\pk{S=1}=c\in (0,1)$ and for some $s <1$ we have $\pk{S<s}=1- c$. Another example is when
\BQN\label{eq:SA}
\pk{S>1-\CC{\frac{1}{u}}}=(1+o(1))c\CC{u^{-\gamma}}, \quad u\to \IF
\EQN
for some $c>0$. In this particular case, the constants $a_n,b_n$ in \eqref{anbbn} can be calculated explicitly as
\BQN\label{add eq2.1}
a_n=(2 \ln   n)^{-\frac{1}{2}}, \qquad b_n=b_{n,\gamma}=(2 \ln n)^{\frac{1}{2}}+(2 \ln
n)^{-\frac{1}{2}}\left( \ln   \varpi-\frac{2\gamma+1}{2}( \ln   \ln n+ \ln   2)\right),
\EQN
with $\varpi=c(2\pi)^{-\frac{1}{2}}\Gamma(1+\gamma)$.

In numerous contributions \CC{(see e.g., Cheng et al. (1998),
 Fahrner and Stadtm\"{u}ller (1998), Cs\'{a}ki and Gonchigdanzan (2002),
Peng et al. (2010), Tan  and Wang (2012), Weng et al. (2012), Hashorva and Weng (2013))} the convergence in distribution for the maxima is strengthen to almost sure convergence. We present such an extension of our main result in the next theorem:

\BT\label{th2.2} Under the \CC{assumptions and} notation of \netheo{th2.1}, if further
\BQN
\rho(n)( \ln   n)^{1+\Delta}( \ln   \ln   n)^{1+\epsilon}=O(1), \quad n\to \IF ,
\label{eq:rho}
\EQN
for some  $\Delta> 2(\gamma-\tau)$ and $\epsilon$ positive, then
for $x \in \R$
\BQN\label{eq2.3}
\frac{1}{ \ln n}\sum^n_{k=1}\frac{1}{k}\I\left(M_k^*\leq a_k x+ b_k \right) \to \Lambda(x), \quad n\to \IF
\EQN
\CC{holds almost surely}, \dE{with $\I(\cdot)$ the indicator function.}
\ET
\BRM i) If \dw{$S$ satisfies \eqref{eq:SA}}, then
Theorem \ref{th2.1} holds under the Berman condition, i.e., we just need to assume therein that  \eqref{eq2.1} is true when  $\Delta=0$.
\dw{Crucial for the proof is that \eqref{eq:bound SA} holds with $\epsilon=0$ if \eqref{eq:SA} holds.}\\
ii) If \eqref{eq:SA} is satisfied and \eqref{eq:rho} holds with
$\Delta=0$, then we have \eqref{eq2.3} also holds with $a_n$ and $b_n$ satisfying \eqref{add eq2.1}.\\
iii) Extension of our results to the case that $X_n,n\ge 1$ is a non-stationary Gaussian sequence is possible.
Various results for extremes of non-stationary Gaussian processes are derived by H\"usler and his co-authors, see for more details
Falk et al. (2010).
\ERM

\section{\ddw{Proofs of the Main Results}}
\prooftheo{th0}
\ccc{The} independence of $S$ and $X$ implies  for any $\nu \in (1,\IF)$ and $u>0$
\begin{equation}\label{wi}
\overline \Phi(u \nu ) \pk{ S> 1/\nu} = \pk{X> u \nu}\pk{ S> 1/\nu} \le  \pk{SX> u, S> 1/\nu} \le \pk{SX> u} \le  \pk{X>u},
\end{equation}
where $\Phi$ is the standard Gaussian distribution on $\R$ and $\ddw{\overline \Phi=1- \Phi}$.
Since for any 
$1<\nu^*< \nu$ we have $\lim_{u\to \IF}  \frac{ \overline \Phi(u \nu )}{\overline \Phi(u \nu^* )}=0$, then
$$ \lim_{u\to \IF}  \frac{ \overline \Phi(u \nu )}{\pk{SX> u}}=\lim_{u\to \IF}  \frac{ \overline \Phi(u \nu )}{\overline \Phi(u \nu^* )} \frac{ \overline \Phi(u \nu^* )}{\pk{SX> u}}
\le \frac{1}{\pk{S> 1/\nu^*}} \lim_{u\to \IF}  \frac{ \overline \Phi(u \nu )}{\overline \Phi(u \nu^* )}
=0$$
implying thus for any $s \in (0,1)$
\BQN
\label{ee}\pk{SX> u} &=& \int_0^ s \overline \Phi(u/x) \, \d F(x) +
\int_s^1 \overline \Phi(u/x) \, \d F(x)\notag\\
& =& O\Bigl(\overline \Phi(u/s)\Bigr)+  \int_s^ 1 \overline \Phi(u/x) \, \d F(x)
\sim  \int_s^ 1 \overline \Phi(u/x) \, \d F(x), \quad u\to \IF.
\EQN
Now,  uniformly for $x\in [1/2,1]$ and some fixed $t\inr$
\begin{equation*}\label{eq2}
\frac{ \overline \Phi(u/x+ (t/x^{2})(x/u))}{\overline \Phi(u/x)} \to
\exp(- t/x^{2}), \quad u\to \IF.
\end{equation*}
Consequently, for $u$ large and
any $\ve^*\in (0,1)$ and $x\in (s,1)$
$$
(1- \ve^*) \overline \Phi(u/x)\exp(-t/s^{2})
\le \overline \Phi(u/x+ t/(xu)) \le (1+ \ve^*) \overline \Phi(u/x)\exp(- t) $$
implying thus for all $u$ large and any $
s\in(1/2,1)$
$$ (1-\ve^{*}) \exp(- t/s^{2}) \le  \frac{\int_s^1 \overline \Phi((u+ t/u)/x) \, \d F(x)}
{\int_s^1 \overline \Phi(u/x) \, \d F(x)} \le (1+ \ve^*)\exp(-t) .$$
Hence for any $\ve\in (0,1)$, since $s$ can be close enough to 1 and
by \eqref{ee}, we obtain
$$ (1-\ve) \exp(- t) \le  \frac{\pk{SX> u+ t/u}}{\pk{SX > u}} \le
(1+ \ve) \exp(- t)$$
and thus $SX$ has distribution function in the Gumbel MDA with scaling function $w(u)=u$.\\
Let $b(t)= \CC{G}^{-1}(1- 1/t)$ with $G^{-1}$ the generalised inverse of the distribution function of $SX$.
In view of \eqref{wi}
for all $t$ large \eH{(write $\Phi^{-1}$ for the inverse of $\Phi$)}
$$ \Phi^{-1}(1- 1/t) \ge b(t)\ge \frac{1}{\nu}\Phi^{-1}\Bigl(1- \frac{1}{t \pk{S> 1/\nu}}\Bigr)$$
\CC{and} since $\nu>1$ can be close enough to 1
$$\CC{b(n) \sim \Phi^{-1}(1- 1/n) \sim (2 \ln n)^{\frac{1}{2}}, \quad a_n \sim (2 \ln n)^{-\frac{1}{2}}, \quad n \to \IF,}$$
hence the claim follows. \QED

\BEL \label{burdon} \dE{Suppose that}  Assumption A holds for $S,S_{\gamma},S_{\tau}$ which are
independent of the random variable $X$ with distribution function \CC{$H$}. \dE{If $H$ has an infinite upper endpoint} and further $H \in GMDA(w)$, then
\BQN\label{eq:bound}
\cZ{\pk{S_{\tau} X> u}\geq \pk{SX> u} \ge  \pk{S_{\gamma} X> u}}
\EQN
holds for all $u$ large\ccc{.} 
\EEL

\prooflem{burdon} By the independence of $S$ and $X$ and the fact that $S$ has \bE{distribution function with} upper endpoint equal 1 for any $\nu>1, u>0$ \bE{we have}
 \BQNY
  \pk{SX> u}= \int_u^{u\CCC{\nu} } \pk{S> u/x} \, \d H(x)+ O(\OH(u\nu )).
  \EQNY
Hence by \eqref{eq:bound S}, for all $u$ large
  $$ \cW{\int_u^{u\nu} \pk{S_{\tau}> u/x} \, \d H(x)}+ O(\OH(u\nu)) \ge \pk{SX> u} \ge  \int_u^{u\nu} \pk{S_{\gamma}> u/x} \, \d H(x)+ O(\OH(u\nu)).$$
A key property of $H  \in GMDA(w)$ is the so-called  {\it Davis-Resnick tail property}, see e.g., Hashorva (2012).
 Specifically, 
 by  Proposition 1.1 of Davis and Resnick (1988)
\BQN\label{Davis}
\lim_{ u \to \IF}  (u w(u))^\mu  \frac{\OH(x u)}{\OH(u)}&=& 0
\EQN
holds for any $\mu\ge 0$ and $x>1$, hence the claim follows now by Theorem 3.1 in Hashorva et al.\ (2010). \QED

\BEL\label{add le3.3}
Let the \cW{positive} random variables $Z_n, n\ge 1$ have  df
$\dw{H}_{n}$ such that \eH{for all large $z$}
\BQN\label{Z_n}
1- \dw{H}_n(z)=\EXP{-\vartheta_n z^q}
\EQN
holds with $q>0,\vartheta_{n}$ positive constants \bE{satisfying
$ \dw{\vartheta}_n \in [a,b], \forall n\ge 1$ with $a<b$ two finite positive constants.}
If further  $Z_n$ is
independent of $S$ which has a regularly varying survival function at 1 with index $\gamma\ge 0$ and
\bE{$u_n,n\ge 1$ are positive constants such that $\limit{n}  u_n=\IF$},
then we have
\BQN\label{eq:SAZ}
\dw{\pk{SZ_n> u_n} \sim}\cW{ \Gamma(\gamma+1)\EXP{-\vartheta_n u_n^q}\pk{S>1-\frac{1}{q\vartheta_n u_n^q}}}.
\EQN
\EEL

\prooflem{add le3.3}
Let $\dw{H}(x)= 1- \exp(- x^q), x>0$ and let $Z$ with distribution function $\dw{H}$ be independent of $S$. By Davis-Resnick tail property of \bE{$\dw{H}$} given in \eqref{Davis} for all large $u_n$, \cW{all $\varepsilon>0$}
\BQNY
\pk{SZ_n>u_n}&=&\int_{u_n}^{\infty}\pk{S>\frac{u_n}{z}}\, \d \dw{H}_{n}(z)\\
&\sim& \int_{u_n}^{u_n(1+\varepsilon)}\pk{S>\frac{u_n}{z}}\, \d \dw{H}_{n}(z)\\
&\sim& \bE{\int_{\vartheta_n^{1/q} u_n}^{\vartheta_n^{1/q}u_n(1+\varepsilon)}\pk{S>\frac{\vartheta_n^{1/q}u_n}{z}}\, \d \dw{H}(z)}\\
&\sim& \pk{SZ>\vartheta_n^{1/q}u_n}\\
&\sim &  \Gamma(\gamma+1)\EXP{-\vartheta_n u_n^q}\pk{S>1-\frac{1}{q\vartheta_n u_n^q}},
\EQNY
where the last step follows from Theorem 3.1 in Hashorva at al. \ (2010). \QED

\BRM
\bE{If}  $S$ \bE{has}  a regularly varying survival function at 1 with index $\gamma\ge 0$,
by the Karamata representation (see e.g., Resnick (1987), p.17),
 we have
\BQNY
\pk{S>1-\frac{1}{q \vartheta_n u_n^q}}\leq c\left(\frac{1}{q \vartheta_n u_n^q}\right)^{\gamma-\epsilon}
\EQNY
with $c>1$, $\epsilon\in (0, \gamma)$. Consequently, by \eqref{eq:SAZ}
\BQN\label{eq:bound SA}
\pk{SZ_n>u_n}\leq c\Gamma(\gamma+1)\EXP{-\vartheta_n u_n^q}\left(\frac{1}{q \vartheta_n u_n^q}\right)^{\gamma-\epsilon}
=O\Bigl( \left(u_n\right)^{-q(\gamma-\epsilon)}\EXP{-\vartheta_n u_n^q}\Bigr)
\EQN
holds for any positive sequence $u_n,n\ge 1$ such that $\limit{n} u_n=\IF$.
\ERM

\BEL\label{add le3.1}
Under the conditions of Theorem \ref{th2.1}, \CC{we have}
\BQN
n\sum_{k=1}^{n-1}|\rho(\CCC{k})|\int_0^{1}\int_0^{1}
\exp\left(-\frac{\left(u_{n}(x)/s\right)^2+\left(u_{n}(x)/t\right)^2}{2(1+|\rho(k)|)}\right)\, \d F(s)\, \d F(t) \to 0,\quad n \to \infty,
\EQN
where
$u_n(x)=a_nx+b_n$ with $a_n$ and $b_n$ are defined in \CC{\eqref{anbbn}} and $x \inr$.
 \EEL

\prooflem{add le3.1}
Denote \eH{$u_{n,c}(z)=a_{n}z+b_{n,c}$}
with $a_{n}$ and $b_{n,c}$ defined in \eqref{add eq2.1}.
\cZ{By Lemma \ref{burdon},
 we have for all large $n$
$b_{n,\gamma}\leq b_n\leq b_{n,\tau}$, hence
$$u_{n,\gamma}(z)\leq u_n(z)\leq u_{n,\tau}(z) \qquad \mbox{and} \qquad \tau\leq \gamma.$$
}
Consequently, \dw{using the assumption of $Z_n$ with $q=2$ and $\vartheta_n=1/(2+2|\rho(n)|), n\ge 1$ in \eqref{Z_n}} and \CCC{\eqref{eq:bound S}}, along the lines of the proof of Lemma 3.2 we obtain for all large $n$
\BQN\label{eq:SY}
\pk{SZ_n>u_n(z)}\leq \pk{SZ_n>u_{n,\gamma}(z)}\leq\pk{S_{\tau}Z_n>u_{n,\gamma}(z)}.
\EQN
\dE{Define next}
$$\sigma=\max_{\ccc{k \ge 1}}|\rho(k)|,  \quad \kappa_n=[n^r],$$
where  $r$ is any positive constant such that
$ r <(1-\sigma)/(1+\sigma)$. This choice of $r$ is possible since by Berman condition and stationarity of the sequence $\sigma < 1$ follows easily. \\
Hereafter $C_1,C_2,C_3$ are positive constants and $\ve \in (0,\tau)$ is taken to be sufficiently small.
By the inequality \eqref{eq:SY} and \eqref{eq:bound SA} (\ddw{denote $F_{\tau}$ the distribution function of $S_{\tau}$}) for all large $n$
\BQNY
&&n\sum_{k=1}^{n-1}|\rho(k)|\int_0^{1}\int_0^{1}
\exp\left(-\frac{\left(u_{n}(x)/s\right)^2+\left(u_{n}(x)/t\right)^2}{2(1+|\rho(k)|)}\right)\, \d F(s)\,\d F(t)\\
&\leq&n\sum_{k=1}^{n-1}|\rho(k)|\int_0^{1}\int_0^{1}
\exp\left(-\frac{\left(u_{n,\gamma}(x)/s\right)^2+\left(u_{n,\gamma}(x)/t\right)^2}{2(1+|\rho(k)|)}\right)
\, \d F_{\tau}(s)\,\d F_{\tau}(t)\\
&\leq&C_1n\sum_{k=1}^{n-1}|\rho(k)|
(u_{n,\gamma}(x))^{-4(\tau-\epsilon)}
\EXP{-\frac{u^2_{n,\gamma}(x)}{1+|\rho(k)|}}\\
&=&C_1n\left(\sum_{k=1}^{\kappa_n}+\sum_{k=\kappa_n+1}^{n-1}\right)|\rho(k)|
(u_{n,\gamma}(x))^{-4(\tau-\epsilon)}\EXP{-\frac{u^2_{n,\gamma}(x)}{1+|\rho(k)|}}=:S_{n1}+ S_{n2}.
\EQNY
According to \CC{\eqref{add eq2.1}} we have
\BQNY
\EXP{-\frac{u_{n,\gamma}^2(x)}{2}}\sim C_2 n^{-1}(u_{n,\gamma}(x))^{1+2\gamma},\quad
u_{n,\gamma}(x)\sim \sqrt{2 \ln   n}, \quad n\to \IF.
\EQNY
As in Lemma 4.3.2 in Leadbetter et al. (1983)
\BQNY
S_{n1}& \le &C_3n^{1+r}(u_{n,\gamma}(x))^{-4(\tau-\epsilon)}
\EXP{-\frac{u_{n,\gamma}^2(x)}{1+\sigma}}\\
&=&O\biggl( n^{1+r}(u_{n,\gamma}(x))^{-4(\tau-\epsilon)}
\left(\frac{(u_{n,\gamma}(x))^{1+2\gamma}}{n}\right)^{\frac{2}{1+\sigma}}\biggr)\\
&=& O\biggl(n^{1+r-\frac{2}{1+\sigma}}( \ln   n)^{\frac{1+2\gamma}{1+\sigma}-2(\tau-\epsilon)}\biggr) \to 0, \quad n\to \IF
\EQNY
\CCC{by} our choice $1+r-\frac{2}{1+\sigma}<0$. Next, with
 $\dE{\sigma}(l)=\max_{\ccc{k \ge l}}|\rho(k)|<1$, we have
\BQNY
S_{n2}
&\leq&C_1 n\dE{\sigma}(\kappa_n)(u_{n,\gamma}(x))^{-4(\tau-\epsilon)}\EXP{-u_{n,\gamma}^2(x)}\sum^{\CCC{n-1}}_{k=\kappa_n+1}
\EXP{\frac{u_{n,\gamma}^2(x)|\rho(k)|}{1+|\rho(k)|}}\\
&\leq&C_1 n^2\dE{\sigma}(\kappa_n)(u_{n,\gamma}(x))^{-4(\tau-\epsilon)}\EXP{-u_{n,\gamma}^2(x)}\EXP{\dE{\sigma}(\kappa_n)u_{n,\gamma}^2(x)}\\
&\CCC{=}&O\biggl(\dE{\sigma}(\kappa_n)(u_{n,\gamma}(x))^{2+4\gamma-4(\tau-\epsilon)}\EXP{\dE{\sigma}(\kappa_n)u_{n,\gamma}^2(x)}\biggr).
\EQNY
By \eqref{eq2.1} \dE{and the fact that $\limit{n} \kappa_n=\IF$} we have for some $\Delta > 2(\gamma-\tau)$
\BQNY
\dE{\sigma}(\kappa_n)(u_{n,\gamma}(x))^{2+2 \Delta }\sim \dE{\sigma}(\kappa_n)(2 \ln   n)^{1+\Delta}
\leq\CCC{\left(\frac{2}{r}\right)}^{1+\Delta }\max_{\ccc{k \ge \kappa_n}}|\rho(k)|( \ln \CCC{k})^{1+\Delta} \to 0
\EQNY
and
$$\dE{\sigma}(\kappa_n)(u_{n,\gamma}(x))^2 \sim 2\dE{\sigma}(\kappa_n) \ln   n
\leq \dE{\sigma}(\kappa_n)(2 \ln   n)^{1+\Delta } \to 0$$
as $n \to \infty$. Since  the exponential term above tends to one and the remaining product
tends to zero, the claim follows. \QED

\prooftheo{th2.1} \CC{Let $\hat{X}_n,n\ge 1$ be independent random variables with the same distribution as $X_1$ and define $\hat{\mathbb{M}}_n^*=\max_{1 \le i \le n}S_i\hat{X}_i$.} If \eqref{eq:bound S} holds, by the independence of the
scaling factors with the Gaussian random variables and Berman's Normal
Comparison Lemma (see e.g., Piterbarg (1996)),
\ddw{and} using Lemma \ref{add le3.1} we
obtain
\begin{eqnarray*}
\lefteqn{\left|\pk{M_n^*\leq u_n(x)}-
\pk{\hat{\mathbb{M}}_n^*\leq u_n(x)}\right|}\\
&\leq&\int_{[0,1]^n}\left|\pk{\bigcap_{k=1}^n
\left\{X_{k}\leq\frac{u_n(z)}{s_k}\right\}}
-\pk{\bigcap_{k=1}^n
\left\{\hat{X}_{k}\leq\frac{u_n(z)}{s_k}\right\}}\right|\, \d F(s_1)\cdots \d F(s_n)\\
&\leq& \frac{1}{4}n\sum_{k=1}^{n-1} |\rho(k)| \int_0^1\int_0^1\EXP{-\frac{(u_n(x)/s)^2+(u_n(s)/t)^2}{2(1+|\rho(k)|)}}\, \d F(s) \d F(t)\\
&\to & 0
\end{eqnarray*}
\cH{as $n\to\infty$}, and thus by Theorem \ref{th0} the claim follows. \QED

\prooftheo{th2.2} In order to show the claim, by Theorem \ref{th2.1} it suffices to prove that
\begin{eqnarray}\label{eq3.6}
\frac{1}{ \ln   n}\sum^n_{k=1}\frac{1}{k}\Bigl(\I\left(M_k^*\leq u_k(x)\right)
-\pk{M_k^*\leq u_k(x)}\Bigr) \to  0, \quad n\to \IF
\end{eqnarray}
holds almost surely, which according to Lemma 3.1 in Cs\'{a}ki and Gonchigdanzan (2002) follows if
for some $\epsilon>0$
\begin{eqnarray*}
\sum^n_{1\leq k< l\leq n}\frac{1}{kl}\Cov\left(\I\left(M_k^*\leq u_k(x)\right),
\I\left(M_l^*\leq u_l(x)\right)\right) = O\biggl( ( \ln   n)^2( \ln   \ln   n)^{-1-\epsilon}\biggr).
\end{eqnarray*}
Next, for any $k<l$ (write below $M_{l,k}^*= \max_{ k< i \le l} S_i X_i$ \CCC{and $\hat{\mathbb{M}}_{l,k}^*=\max_{ k< i \le l} S_i \hat{X}_i$})
\begin{eqnarray*}
&&\Cov\left(\I\left(M_k^*\leq u_k(x)\right),
\I\left(M_l^*\leq u_l(x)\right)\right)\\
&\leq&2\E\left|\I\left(M_l^*\leq u_l(x)\right)-\I\left(M_{l,k}^*\leq u_l(x)\right)\right|
+\left|\Cov\left(\I(M_k^*\leq u_k(x)),\I(M_{l,k}^*\leq u_l(x))\right)\right|\\
&\leq&2\left|\pk{M_{l,k}^*\leq u_l(x)}
-\pk{\hat{\mathbb{M}}_{l,k}^*\leq u_l(x)}\right|
+2\left|\pk{M_{l}^*\leq u_l(x)}
-\pk{\hat{\mathbb{M}}_{l}^*\leq u_l(x)}\right|\\
&&+2\left|\pk{\hat{\mathbb{M}}_{l,k}^*\leq u_l(x)}
-\pk{\hat{\mathbb{M}}_{l}^*\leq u_l(x)}\right|\\
&&+\left|\pk{M_k^*\leq u_k(x),M_{l,k}^*\leq u_l(x)}
-\pk{M_k^*\leq u_k(x)}\pk{M_{l,k}^*\leq u_l(x)}\right|\\
&=&P_1+P_2+P_3+P_4.
\end{eqnarray*}
In view of Berman's Normal Comparison Lemma \CC{and \eqref{eq:rho}, along the same lines of the proof of Lemma \ref{add le3.1}}, we have
$$P_i= O \biggl(  ( \ln   \ln n)^{-1-\epsilon}\biggr),  \quad  \ i=1,2,4.$$
Further, since
\BQNY
P_3=\mathbb{G}^{l-k}(u_l(x))-\mathbb{G}^{l}(u_l(x))\leq \frac{k}{l},
\EQNY
where $\mathbb{G}$ is the df of $S_1\hat{X}_1$, we establish the claim. \QED

\cH{\textbf{Acknowledgments}}: Z. Weng kindly acknowledges financial support by the Swiss National Science Foundation Grant 200021-134785 and the project RARE -318984, a Marie Curie FP7 IRSES Fellowship.
E. Hashorva kindly acknowledges partial support by the Swiss National Science Foundation Grant 200021-140633/1.

\end{document}